\documentclass[graybox]{svmult}

\usepackage{mathptmx}       
\usepackage{helvet}         
\usepackage{courier}        
\usepackage{type1cm}        
\usepackage{makeidx}         
\usepackage{graphicx}        
\usepackage{multicol}        
\usepackage[bottom]{footmisc}

\usepackage{amssymb}
\usepackage{amsmath}

\makeindex             

\begin{document}
\title*{Iteration scheme for initial value problem for PDEs: Existence, convergence and comparison}
\titlerunning{Iteration scheme for PDEs: Existence, convergence and comparison}
\author{Josef Rebenda and Zden\v{e}k \v{S}marda}
\institute{Josef Rebenda \at CEITEC BUT, Brno University of Technology, Technicka 3058/10, 61600 Brno, Czech Republic, \email{josef.rebenda@ceitec.vutbr.cz}
\and Zden\v{e}k \v{S}marda \at CEITEC BUT, Brno University of Technology, Technicka 3058/10, 61600 Brno, Czech Republic, \email{smarda@feec.vutbr.cz}}
%
%
\maketitle

\abstract*{Each chapter should be preceded by an abstract (10--15 lines long) that summarizes the content. The abstract will appear \textit{online} at \url{www.SpringerLink.com} and be available with unrestricted access. This allows unregistered users to read the abstract as a teaser for the complete chapter. As a general rule the abstracts will not appear in the printed version of your book unless it is the style of your particular book or that of the series to which your book belongs.
Please use the 'starred' version of the new Springer \texttt{abstract} command for typesetting the text of the online abstracts (cf. source file of this chapter template \texttt{abstract}) and include them with the source files of your manuscript. Use the plain \texttt{abstract} command if the abstract is also to appear in the printed version of the book.}

\abstract{Results about existence and uniqueness of solutions of initial value problem for certain types of partial differential equations are recalled as well as iterative scheme and an error estimate for approximate solutions obtained using this scheme. Several numerical examples are presented to demonstrate how the proposed iterative scheme can be applied, with emphasis given to verifying assumptions of using the scheme. Comparison to other recently presented results is done in this respect.}

\section{Introduction}
\label{sec:1}

The significance of partial differential equations is growing in the field of mathematical analysis recently.
Researchers are looking for new methods which are computationally efficient and simple in application to find exact or approximate solutions of partial differential equations.

Many techniques
including variational iteration method, homotopy analysis method, homotopy perturbation method and differential transformation method have been recently used for solving particular problems. Unfortunately, most of the above mentioned methods need some analytical preparation before using some computer software. Most of the algorithms require software which is able to perform symbolical integration.

Following direction of current research, we found a lot of papers whose content was unsatisfactory. Some authors do not verify conditions that justify using of the particular method. Some authors even do not mention any assumptions under which the method is applicable. Some also try to prove existence and uniqueness of solutions of differential equations in a Banach space of functions which are only continuous on some compact set.

Confused by ongoing research in this area, we started to examine this approach in \cite{rebenda:1}. However, we found out that the conditions of the approach are too restrictive. Therefore, we introduced milder and more reasonable sufficient conditions  in more general situation of initial value problems for multidimensional PDEs in paper \cite{rebenda:2}, where we proved local existence and uniqueness of solution of considered problem using Banach fixed-point theorem. An iterative scheme for solving initial problems for certain types of PDEs was derived in the paper as well as several examples.

The main purpose of this paper is to supplement the paper \cite{rebenda:2} with comparison to other recently published results dealing with methods producing iteration algorithms.

The paper is organized as follows. In Section \ref{sec:2}, we repeat all necessary notations, definitions and conditions which are introduced in paper \cite{rebenda:2}. In view of this preparatory part we recall the main results as well as corollaries in Section \ref{sec:3}. Being novel and the most important part of this paper, Section \ref{sec:4} consists of comparison to other results and several explanatory examples. Conclusions are summarized in Section \ref{sec:5}.

\section{Preliminaries}
\label{sec:2}
Let $\Omega$ be a compact subset of $ \mathbb{R}^k$. Denote $J = [-\delta,\delta] \times \Omega$, where $\delta>0$ will be specified later, then $J$ is a compact subset of $ \mathbb{R}^{k+1}$. Let $u(t,x)=u(t,x_1,\dots,x_k)$
be a real function of $k+1$ variables defined on $J$. We introduce the following operators: $\nabla = (\frac{\partial}{\partial t},\frac{\partial}{\partial x_1}, \dots,
\frac{\partial}{\partial x_k})$ and $D = (\frac{\partial}{\partial x_1}, \dots, \frac{\partial}{\partial x_k})$. We deal with partial differential equations
\begin{equation}\label{eq1}
\frac{\partial^n}{\partial t^n} u(t,x) = F (t, x, u, \nabla u, \dots, \nabla^m u) \quad \text{for } m<n
\end{equation}
and
\begin{equation}\label{eq2}
\frac{\partial^n}{\partial t^n} u(t,x) = F (t, x, u, \nabla u, \dots, \nabla^{n-1} u, D \nabla^{n-1} u, D^2 \nabla^{n-1} u, \dots, D^{m-(n-1)} \nabla^{n-1} u),  m \geq n.
\end{equation}
In both cases, left-hand side of the equation contains only the highest derivative with respect to $t$. We do not consider equations where the order of partial derivatives
with respect to $t$ is $n$ or higher on the righthand side, including mixed derivatives.

When convenient, we will use multiindex notation as well:
$$
\frac{\partial^{\vert \alpha \vert}}{\partial x^{\alpha}} = \frac{\partial^{\vert \alpha \vert}}{\partial x_1^{\alpha_1} \partial x_2^{\alpha_2} \dots \partial x_k^{\alpha_k}},
$$
where $\vert \alpha \vert = \alpha_1 + \alpha_2 + \dots + \alpha_k$.

Denote $N=\max \{ m,n \}$.
We consider equation \eqref{eq1} or \eqref{eq2} with the set of initial conditions
\begin{align}\label{eq3}
u(0,x) &= c_1 (x),\notag \\
\frac{\partial}{\partial t} u(0,x) &= c_2 (x), \notag \\
 &\vdots \notag \\
\frac{\partial^{n-1}}{\partial t^{n-1}} u(0,x) &= c_n (x),
\end{align}
where initial functions $c_i (x)$, $i=1,\dots,n$ are taken from space $C^N (\Omega,  \mathbb{R})$. It means that we are looking for classical solutions.

For the purpose of clarity, we emphasize that our formulation covers for instance heat, wave, Burger, Boussinesq or Korteweg-de Vries (KdV) equations.

Obviously, $F \colon J \times  \mathbb{R}^k \times  \mathbb{R} \times  \mathbb{R}^{k+1} \times  \mathbb{R}^{(k+1)^2} \times \dots \times  \mathbb{R}^{(k+1)^m} \to  \mathbb{R}$ if $m<n$,
and $F \colon J \times  \mathbb{R}^k \times  \mathbb{R} \times  \mathbb{R}^{k+1} 
\times  \mathbb{R}^{(k+1)^2} \times \dots \times  \mathbb{R}^{(k+1)^{n-1}} \times  \mathbb{R}^{k(k+1)^{n-1}} \times  \mathbb{R}^{k^2(k+1)^{n-1}} \times \dots \times  \mathbb{R}^{k^{m-n+1}(k+1)^{n-1}} \to  \mathbb{R}$ if $m \geq n$.
Denote 
\begin{equation}
K_1 = \frac{(k+1)^{m+1} -1}{k} \quad \text{for } m<n
\end{equation}
and
\begin{equation}
K_2 = \frac{(k+1)^{n-1} -1}{k} + (k+1)^{n-1} \frac{k^{m-n+2} -1}{k-1} \quad \text{for } m \geq n.
\end{equation}
Then, if we consider $u$ as dependent variable, we see that $F$ is a function of $k+1+K_1$ variables in case $m<n$ or $k+1+K_2$ variables in case $m \geq n$.

Denote
\begin{equation}\label{eq4}
u_0 (t,x) = \sum_{i=1}^n c_i (x) \frac{t^{i-1}}{(i-1)!} = \sum_{i=1}^n \left( \frac{\partial^{i-1}}{\partial t^{i-1}} u(0,x) \right) \frac{t^{i-1}}{(i-1)!}.
\end{equation}
Then $u_0 \in C^N (J,  \mathbb{R})$.

We suppose that $F$ is Lipschitz continuous in last $K_1$ ($m<n$) or $K_2$ ($m \geq n$) variables, i.e. $F$ satisfies condition
\begin{equation}\label{eq5}
\vert F (t,x,y_1, \dots, y_{K_l}) - F (t,x,z_1, \dots, z_{K_l}) \vert \leq L \left( \sum_{i=1}^{K_l} \vert y_i - z_i \vert \right), \quad l = 1 \text{ or } 2,
\end{equation}
on a compact set which is defined as follows: There is $R \in  \mathbb{R}, R>0$ such that \eqref{eq5} holds on 
\begin{equation}\label{eq6}
J \times \prod_{\alpha_0 + \vert \alpha \vert \leq m}
[c_{\alpha_0, \alpha}, d_{\alpha_0, \alpha} ],
\end{equation}
where
\begin{equation}\label{eq7}
c_{\alpha_0, \alpha} = \min_{(t,x) \in J} \left[ \frac{\partial^{\alpha_0 + \vert \alpha \vert}}{\partial t^{\alpha_0} \partial x^\alpha} u_0 (t,x) \right] -R,
\quad
d_{\alpha_0, \alpha} = \max_{(t,x) \in J} \left[ \frac{\partial^{\alpha_0 + \vert \alpha \vert}}{\partial t^{\alpha_0} \partial x^\alpha} u_0 (t,x) \right] +R,
\end{equation}
and $\alpha_0 <n$ in all cases.

Since $F$ is continuous on compact set, $\vert F \vert$ attains its maximal value on this set, denote it $M$. Then we put
\begin{equation}\label{eq10}
\delta = \left( \frac{R \cdot (n-1)!}{M} \right)^{1/n}.
\end{equation}

Finally, define the following operator
\begin{equation}\label{eq8}
T u(t,x) = u_0 (t,x) + \int_0^t \frac{(t-\xi)^{n-1}}{(n-1)!} F \bigl(\xi, x, u (\xi,x), \nabla u(\xi,x), \dots \bigr) d \xi,
\end{equation}
where the function $F$ has either $k+1+K_1$ or $k+1+K_2$ arguments and the last $K_1$ or $K_2$ arguments involve dependent variable $u$.

\section{Overview of results}
\label{sec:3}
This section is devoted to recalling the results presented in \cite{rebenda:2}. We would like to emphasize that all results have local character.

\begin{theorem}\label{th1}
Let the condition \eqref{eq5} hold. Then problem consisting of equation \eqref{eq1} or \eqref{eq2} and initial conditions \eqref{eq3} has a unique local solution on $(-\delta, \delta) \times \Omega$, where $\delta$ is defined by \eqref{eq10}.
\end{theorem}

\begin{theorem}\label{th2}
Assume that condition \eqref{eq5} holds. Then iterative scheme $u_p = T u_{p-1}$, $p \geq 1$ with initial approximation $u_0$ defined by \eqref{eq4}, where $T$ is defined by \eqref{eq8}, converges to unique local solution $u(t,x)$ of problem \eqref{eq1}, \eqref{eq3}, respective \eqref{eq2}, \eqref{eq3}. Moreover, we have the following error estimate for this scheme:
\begin{equation}\label{eq13}
\Vert u(t,x) - u_p (t,x)\Vert_{C^N} \leq \frac{R \cdot \gamma^p}{1-\gamma}
\end{equation}
on $(-\delta_1, \delta_1) \times \Omega$, where $\delta_1$ is chosen such that operator $T$ is a contraction, $\gamma = \frac{L \cdot \delta_1^n}{(n-1)!}$ and constants $L$ and $R$ are defined by \eqref{eq5} and \eqref{eq6}.
\end{theorem}

\begin{corollary}\label{c1}
Let condition \eqref{eq5} be valid and suppose that $F$ can be written as $G+g$:
\begin{equation}\label{eq12}
F \bigl(t, x, z (t,x), \nabla z(t,x), \dots \bigr) = G \bigl(t, x, z (t,x), \nabla z(t,x), \dots \bigr) + g(t,x).
\end{equation}
Then we may choose initial approximation
\begin{align}\label{eq11}
\bar{u}_0 &= u_0 + \int_0^t \frac{(t-\xi)^{n-1}}{(n-1)!} g(\xi, x) d \xi\\
&= \sum_{i=1}^n \left( \frac{\partial^{i-1}}{\partial t^{i-1}} u(0,x) \right) \frac{t^{i-1}}{(i-1)!} + \int_0^t \frac{(t-\xi)^{n-1}}{(n-1)!} g(\xi, x) d \xi.
\end{align}
\end{corollary}

\section{Comparison}
\label{sec:4}

Now we are prepared to compare our results to other recently published results. We would like to point out that formulation and verification of conditions sufficient to use particular algorithms is important.

\begin{example}\label{ex1}
The first example to compare was published in \cite{sweilam:1} in 2010 as Example 1. Consider the following third-order nonlinear PDE
\begin{equation}\label{eqex1}
u_t + 6 u^2 u_x + u_{xxx} = 0
\end{equation}
with initial function $u(0,x) = \sqrt{c} \ {\rm sech} (k+\sqrt{c} x)$, where $c \geq 0$ and $k$ is a constant.

The authors of paper \cite{sweilam:1} claim that their "variational iteration formula" converges "for some constants $\gamma_i = \alpha_i + m_i T < 1$". However, no verification of validity of such condition is done in the examples.

On the other hand, to apply our iteration scheme, we need to verify that condition \eqref{eq5} is valid on some connected compact set of the form \eqref{eq6}. In our terms, $F(t,x,y_1,y_2,y_3) = -6 y_1^2 y_2- y_3$. Obviously, this function has continuous partial derivatives of all orders in $\mathbb{R}^5$, hence the difference $|F(t,x,y_1,y_2,y_3) - F(t,x,z_1,z_2,z_3)|$ can be estimated by $M_1 |y_1 - z_1| + M_2 |y_2 - z_2| + M_3 |y_3 - z_3|$, where $M_i$, $i=1,2,3$ are maxima of partial derivatives of function $F$ with respect to last three variables on \eqref{eq6}. It means that condition \eqref{eq5} is fulfilled and our iteration scheme can be applied to find local approximation of unique solution of considered problem.
\end{example}

\begin{example}\label{ex2}
The second example is taken from paper \cite{odibat:1} published in 2010 where it is listed as Example 2. Consider a first-order nonlinear ODE (which we may regard as PDE)
\begin{equation}\label{eqex2}
u'(t) + u^3 (t) = t^3 +3t^2 + 3t +2
\end{equation}
subject to initial condition $u(0) = 1$.

In this paper, the author does the verification of sufficient condition. However, this verification leads to a false conclusion that the algorithm does not converge.

In our terms, $F(t,y_1) = t^3 +3t^2 + 3t +2 - y_1^3$. Again, this function has continuous partial derivatives of all orders in $\mathbb{R}^2$, hence the same argument as in Example \ref{ex1} shows that our iteration scheme (which in fact is Picard's successive approximations scheme in this case) is applicable in a neighbourhood of origin. Indeed, calculating first few terms, we get the following sequence:
\begin{eqnarray*}
u_1 & = & 1 + t + \frac{3}{2} t^2 + t^3 + \frac{1}{4} t^4,\\
u_2 & = & 1 + t \quad \quad - \frac{3}{2} t^3 - 3 t^4 + \ldots ,\\
u_3 & = & 1 + t \quad \quad \quad \quad + \frac{9}{8} t^4 + \ldots .
\end{eqnarray*}
Calculating a few more terms reveals that all terms of order $2$ and higher are vanishing one by one in each step, thus the sequence converges to exact solution $u(t) = 1 +t$ near the origin. This also is the case in \cite{odibat:1}, if the author tried to calculate the sequence of partial sums. However, the verification led to the false conclusion which prevented him to try it.
\end{example}

\begin{example}\label{ex3}
The third example was published in \cite{he:1} in 1997 as Example 4. It is a first order nonlinear PDE
\begin{equation}\label{eqex3}
u_t = x^2 - \frac{1}{4} u_x^2
\end{equation}
with initial function $u(0,x) = 0$.

In this case, since the cited paper does not contain any assumptions, conditions, theorems or proofs, no verification was done in the examples.

The same argument as in previous two examples (continuous partial derivatives of $F$) allows us to apply our iteration scheme to find local approximation of unique solution of the problem. The first few approximations are
\begin{eqnarray*}
u_1 (t,x) & = & x^2 t,\\
u_2 (t,x) & = & x^2 (t - \frac{1}{3} t^3),\\
u_3 (t,x) & = & x^2 (t - \frac{1}{3} t^3 + \frac{2}{15} t^5 + o(t^5)),\\
u_4 (t,x) & = & x^2 (t - \frac{1}{3} t^3 + \frac{2}{15} t^5 + \frac{17}{315} t^7 + o(t^7)).
\end{eqnarray*}
Careful examination of this sequence leads to a guess that it converges to the function $u(t,x) = x^2 \ {\rm tanh} \ t$. Indeed, it is not difficult to verify that this function is an exact solution of the considered problem. According to Theorem \ref{th1}, this solution is unique, hence our sequence of approximations indeed converges to this unique solution.
\end{example}

\section{Conclusion}
\label{sec:5}
We formulated an existence and uniqueness result for initial value problem for certain classes of PDEs in this paper. Using this result, we derived an iteration scheme, investigated its convergence and compared our results to other recently published results. It turned out that our results are less restrictive than the compared results. A subject of further investigation is to develop similar approach for systems of PDEs and for other types of problems.

\begin{acknowledgement}
This research was carried out under the project CEITEC 2020 (LQ1601) with financial support from the Ministry of Education, Youth and Sports of the Czech Republic under the National Sustainability Programme II.
This work was also supported by the Czech Science Foundation under the project 16-08549S. The work was realized in CEITEC - Central European Institute of Technology with research infrastructure supported by the project CZ.1.05/1.1.00/02.0068 financed from European Regional Development Fund.
This support is gratefully acknowledged. 
\end{acknowledgement}

%
%
%

\end{document}